\documentclass[11pt,a4paper,reqno]{amsart}

\usepackage{mathpazo}
\usepackage{mdframed}
\usepackage[shortlabels]{enumitem}
\usepackage{tikz-cd}
\usepackage{appendix}
\usepackage{graphicx}
\usepackage{caption}
\usepackage{amsthm,amsmath,tikz}
\usepackage{mathrsfs}
\usepackage{xspace}
\usepackage[stable,perpage]{footmisc}
\usepackage{marginnote}
\usepackage{amsfonts,amsmath,amssymb}
\usepackage{fancyhdr}
\usepackage{tikz-cd}
\usetikzlibrary{patterns,decorations.pathmorphing,shapes}
\usepackage{mathtools,stmaryrd}
\usetikzlibrary{decorations}
\usepackage{pgfplots}
\usepackage{hyperref,setspace,microtype}
\usepackage[margin=1in]{geometry}

\tikzset{
  symbol/.style={
    draw=none,
    every to/.append style={
      edge node={node [sloped, allow upside down, auto=false]{$#1$}}}
  }
}

\usepackage{mdframed}
\usepackage{setspace}

\theoremstyle{plain}
\newtheorem{thm}{Theorem}[section]
\newtheorem{lem}[thm]{Lemma}
\newtheorem{prop}[thm]{Proposition}
\newtheorem{cor}[thm]{Corollary}

\theoremstyle{definition}
\newtheorem{defn}[thm]{Definition}

\newtheorem{ex}[thm]{Example}

\theoremstyle{remark}
\newtheorem*{pf}{Proof}
\newtheorem{rem}[thm]{Remark}
\newtheorem*{note}{Note}
\newcommand{\N}{\mathbb{N}}

\newcommand{\R}{\mathbb{R}}

\renewcommand{\O}{\mathcal{O}}
\renewcommand{\P}{\mathbb{P}}

\newcommand{\A}{\mathbb{A}}
\newcommand{\G}{\mathbb{G}}

\newcommand{\Spec}{\text{Spec}}

\newcommand{\etale}{\text{\'etale}\xspace}

\newcommand{\ACGS}{\text{ACGS}}

\newcommand{\vir}{\text{vir}}
\newcommand{\Glog}{\mathbb{G}_\text{log}}
\newcommand{\Mgp}{M^\text{gp}}
\newcommand{\VZ}{\mathcal{V\mkern-3mu Z}}

\newcommand{\rad}{\text{rad}}

\usepackage{stackengine}
\stackMath
\newcommand\tsup[2][2]{%
 \def\useanchorwidth{T}%
  \ifnum#1>1%
    \stackon[-.5pt]{\tsup[\numexpr#1-1\relax]{#2}}{\scriptscriptstyle\sim}%
  \else%
    \stackon[.5pt]{#2}{\scriptscriptstyle\sim}%
  \fi%
}

\newcommand{\sqC}{\scalebox{0.8}[1.2]{$\sqsubset$}}

\usepackage{afterpage}

\title{Moduli of elliptic curves in products of projective spaces}
\author{Wanlong Zheng}
\address{Department of Pure Mathematics and Mathematical Statistics, University of Cambridge}
\email{wz302@cam.ac.uk}
\date{\today}

\begin{document}

\maketitle
\begin{abstract}
   We exhibit a smooth compactification of the moduli space of elliptic curves in a product of projective spaces with tangency along a subset of its toric boundary divisors. This is a Vakil--Zinger type of desingularization for maps to a product of projective spaces using ideas of elliptic singularities and logarithmic geometry, extending the recent work by Ranganathan--Santos-Parker--Wise. We use this to construct the virtual fundamental classes of the spaces of genus 1 maps to a special class of simple normal crossings pairs.
\end{abstract}
\setcounter{tocdepth}{1}
\tableofcontents
\section{Introduction}

\subsection{Main results}

Let $n_1,\dots,n_a$ be positive integers, and $\P^{\underline{n}}:=\P^{n_1}\times\cdots\times\P^{n_a}$. Let $D$ be a subset of $\partial \P^{\underline{n}}$, the set of all irreducible toric boundary divisors, and $\Gamma$ be the matrix where $\Gamma_{ij}$ is the contact order (tangency order) of the $i$-th marked point with a component $D_j$ of $D$.

Denote by $\ACGS^\Gamma_{1,n}(\P^{\underline{n}},D)$ the moduli space of stable logarithmic maps in genus 1 with respect to the pair $(\P^{\underline{n}},D)$, constructed by Abramovich--Chen and Gross--Siebert in~\cite{AC, GS, ACGS}.

\begin{thm}\label{thm:big}
    There is a proper moduli space $\VZ_{1,n}^\Gamma(\P^{\underline{n}}, D)$
    of \textbf{well--spaced} stable maps from $n$-pointed genus 1 curves to the pair $(\P^{\underline{n}}, D)$ with contact orders $\Gamma$. It is logarithmically smooth, and the morphism
    $$\VZ_{1,n}^\Gamma(\P^{\underline{n}}, D)\to \ACGS^\Gamma_{1,n}(\P^{\underline{n}}, D)$$ maps the moduli space birationally onto the main component of $\ACGS^\Gamma_{1,n}(\P^{\underline{n}},D)$.
\end{thm}

The well-spacedness condition removes certain strata that are of excess dimension, yielding a \textit{reduced} theory in genus 1. It preserves the main component of the usual moduli space of stable maps to the target pair. We will use the term \textit{logarithmic desingularization} to indicate that the resulting space $\VZ_{1,n}^\Gamma(\P^{\underline{n}}, D)$ is logarithmically smooth.

By taking $D=\emptyset$ in Theorem \ref{thm:big}, the map becomes $$\VZ_{1,n}(\P^{\underline{n}}, \emptyset) \to \overline{\mathcal{M}}_{1,n}(\P^{\underline{n}},\boldsymbol{d}),$$ which is a Vakil--Zinger type of logarithmic desingularization~\cite{VZ} of the main component of the moduli space of stable logarithmic maps to a product of projective spaces in genus 1. Although the techniques in~\cite{VZ} may desingularize this space directly, a detailed construction has not yet appeared in the literature to our knowledge. In particular, we show in Section \ref{sec:vzexample} that the Vakil--Zinger construction does not desingularize the moduli space of $(2,2)$-curves in $\P^1\times \P^1$. Furthermore, our construction is a modular compactification, whereas the Vakil--Zinger desingularization does not have a modular interpretation.

We extend the construction to more general target spaces. Let $X$ be a smooth projective variety, and $Y=\sum Y_i$ a simple normal crossings divisor whose irreducible components $Y_i$ are defined by sections of the same line bundle $\phi^*\mathcal{O}(1)$ for some embedding $\phi$ and hyperplanes $H_i$:
$$\phi:(X,\sum Y_i)\to (\P^r,\sum H_i).$$
We call this a \textit{very ample SNC pair of the same degree} $(X,Y)$. A typical example is the logarithmic Calabi-Yau pair $X=\P^1\times \P^1$ with $Y=Y_1+Y_2$ a sum of two $(1,1)$-curves.

\begin{thm}\label{thm:1.1}
  There exists a proper moduli stack $\VZ_{1,n}^\Gamma(X,Y)$ of well--spaced genus 1 stable maps to the pair $(X,Y)$ with contact orders $\Gamma$. It carries a virtual fundamental class in the expected degree.
\end{thm}

Because of the well-spacedness condition, the virtual class of $\VZ_{1,n}^\Gamma(X,Y)$ isolates the contribution from the main component of the usual space of stable maps.

Lastly, we provide a comparison result of the virtual classes when forgetting a certain component and several markings. This comparison result can be used to recursively compute the virtual classes and invariants of some moduli spaces. It, therefore, shows that the genus 1 theory defined in this paper is self-consistent, where the calculations of the invariants can be achieved without leaving this genus 1 theory.

Let $(X,Y)$ be as above, and $Y_1$ a component of $Y$. Write $Y-Y_1$ for the divisor obtained by forgetting $Y_1$.
\begin{thm}[Consistency]
  For the pair $(X,Y)$ and the component $Y_1$, if both of the following are satisfied:
  \begin{itemize}
    \item $\Gamma_{i1}=1$ for $1\le i\le m$ and $\Gamma_{i1}=0$ for $i>m$; namely only $\{p_1,\dots,p_m\}$ intersect $Y_1$, and they intersect transversely; and
    \item $\Gamma_{ij}=0$ for $1\le i\le m$ and $j>1$; namely those $m$ points do not intersect any other $Y_j$,
  \end{itemize}
  then there is a forgetful map
  $$g:\VZ_{1,n}^\Gamma(X, Y) \to \VZ_{1,n-m}^{\Gamma'}(X, Y-Y_1)$$
  that relates the virtual classes by the identity:
  $$g_*[\VZ_{1,n}^\Gamma(X, Y)]^{\text{vir}}=m!\cdot[\VZ_{1,n-m}^{\Gamma'}(X, Y-Y_1)]^{\text{vir}}$$
  where $\Gamma'$ is the submatrix of contact orders by forgetting the $m$ points and $Y_1$.
\end{thm}
\subsection{Context and overview}
This paper concerns the curve counting in genus 1. The theory developed here combines several viewpoints from elliptic singularities, tropical geometry, logarithmic geometry and the traditional theory of moduli spaces.

Ranganathan--Santos-Parker--Wise recently in~\cite{RSW1} reinterpreted the Vakil--Zinger desingularization of moduli space of genus 1 stable maps to projective spaces using tools developed in logarithmic and tropical geometry. More precisely, they constructed a proper and smooth moduli space $\VZ_{1,n}(\P^n,d)$ of stable maps from radially aligned genus 1 curves to $\P^n$ satisfying the factorization condition, and the map
$$\VZ_{1,n}(\P^n,d)\to \overline{\mathcal{M}}_{1,n}(\P^n,d)$$
to the space of genus 1 stable maps desingularizes the main component of the target space. Moreover, the space $\VZ_{1,n}(\P^n,d)$ is birational to the Vakil--Zinger blowup construction, and it has the advantage of being a moduli problem (or indeed, using a slightly more refined notion of \textit{centrally aligned} curves will give an isomorphism; however we will not be using this notion in this paper).

We consider the \textit{relative} problem of maps to a target pair $(X,Y)$, where $Y$ is a divisor of $X$, with prescribed contact orders $\Gamma$ of the marked points to $Y$. Ranganathan--Santos-Parker--Wise extended their previous results to maps to any toric pair target in the sequel paper~\cite{RSW2}. More precisely, there is a proper and toroidal moduli space $\VZ_{1,n}^\Gamma(X,\partial X)$ of stable maps from genus 1 curves to a toric variety $X$ with contact orders $\Gamma$ satisfying the \textit{well--spacedness condition}. This is a technical condition which we will briefly explain in Section \ref{sec:prelim} and generalize in Section \ref{sec:moduli}. Similar to the absolute case, the map
$$\VZ_{1,n}^\Gamma(X,\partial X)\to \ACGS^\Gamma_{1,n}(\P^{\underline{n}},\partial X)$$
is a logarithmic desingularization of the main component of the space of logarithmic stable maps to the pair $(X,\partial X)$.

The results in this paper further generalize the above situation to any very ample SNC pair of the same degree. In particular, by choosing the pair $(\P^{\underline n},\emptyset)$ and applying Theorem \ref{thm:big}, we have an explicit construction of the Vakil--Zinger type of desingularization for genus 1 stable maps to a product of projective spaces. This is an explicit desingularization construction for the product of projective spaces target.

An intermediate case where the target pair is $(\P^r,H)$ for a smooth hyperplane $H$ is examined by Battistella--Nabijou--Ranganathan in~\cite{BNR}. They constructed a proper moduli space $\VZ_{1,n}^\Gamma(\P^r,H)$ and showed it is toroidal by carefully analyzing the obstructions to stable maps. In this sense, we provide a shorter and more geometric proof to this fact in this paper by reducing it to the smoothness of the spaces constructed in~\cite{RSW2}.

Towards the end of the paper, we show the theory developed here is self-consistent. We identify a case where certain marked points together with a component of the divisor, which are called \textit{fictitious}, do not meaningfully contribute towards the Gromov--Witten invariants calculation. This result is similar to the one proposed by Gathmann in~\cite{gathhyper}. There is an explicit relation on the virtual fundamental classes along the forgetful map, and it enables one to calculate certain Gromov--Witten invariants recursively.

\subsection{Future directions}
An immediate generalization would be to remove the \textit{same degree} requirement for the very ample SNC pair $(X,Y)$. The main problem for applying the same techniques discussed in this paper to this case is the lack of virtual classes on $\VZ_{1,n}(X,\emptyset)$. See Section \ref{sec:generalpair} for some further discussions.

We could also consider curves with higher genus. One central idea we used in genus 1 is to contract an elliptic curve in a family to get an elliptic Gorenstein singularity. In genus 2, Battistella--Carocci in~\cite{BC} were able to classify all genus 2 Gorenstein singularities, and constructed a modular desingularization of $\overline{\mathcal M}_{2,n}(\P^r,d)^{\text{main}}$. The geometry is more delicate, with both isolated and non-reduced singularities appearing in the contraction. More generally, Bozlee in ~\cite{Boz} constructed a contraction of multiple subcurves of any genus inside of a family of any genus. The singularities that can appear include the elliptic Gorenstein singularities we considered, but there are others.

The relative geometry is necessarily more delicate in higher genus, but we hope our results in genus 1 may provide some hints to understanding the main component of the toric contact cycles~\cite{samdhruv} and the moduli spaces themselves.

Finally, although we have provided a consistency result that makes certain calculations possible, we still expect the more powerful tools such as the degeneration formula and the toric localization to hold in the proposed genus 1 theory.

\subsection*{Acknowledgments} I would like to thank my advisor Dhruv Ranganathan for introducing and suggesting this topic, and for his continuous support and guidance. I would also like to thank Navid Nabijou for explaining various concepts and some very helpful discussions. Finally I learned a lot from talking to various people from the research group at Cambridge, especially with Patrick Kennedy-Hunt and Qassim Shafi.

The author is supported by Cambridge International Trust and DPMMS at the University of Cambridge.

\section{Preliminaries}\label{sec:prelim}
\subsection{Logarithmic curves and tropicalizations}
We review the logarithmic structure on a curve and the construction of its associated tropicalization. A more careful and detailed approach can be found in~\cite{Kato}.

A logarithmic scheme $(S,M_S)$ is a scheme $S$ with a sheaf of monoids $M_S$ and a map of sheaves of monoids
  $$\alpha:M_S\to \mathcal{O}_X$$
such that $\alpha$ restricts to a bijection $\alpha^{-1}(\mathcal{O}_X^\times)\cong \mathcal{O}_X^\times$.

The \textit{ghost sheaf} of $S$ is defined by $\overline{M}_S:=M_S/\mathcal{O}_X^\times$. Similarly define $\overline{M}_S^{\text gp}:=M_S^{\text gp}/\mathcal{O}_X^\times$.

\begin{defn}
  A \textbf{family of logarithmic curves} over a logarithmic scheme $(S,M_S)$ is a logarithmically smooth and proper morphism
  $$\pi:(C,M_C)\to (S,M_S)$$
  of logarithmic schemes, such that the fibers are 1--dimensional and connected, and $\pi$ is integral and saturated.
\end{defn}

We often call this a \textit{logarithmic curve over} $S$ for short, suppressing the sheaves of monoids.

Associated to a logarithmic curve $C\to S$ is a family of tropical curves. This is achieved in two steps.

\begin{itemize}
  \item If $S=\Spec(P\to k)$ is a logarithmic point associated to a toric monoid $P$, then the logarithmic structure keeps track of a generalized edge length $l_e\in P=\overline{M}_S$. Together with the genera of each component of $C$, the data can be packaged into its \textit{tropicalization} as a map from a cone complex $\sqC$ to the dual cone $\sigma_P$
  \[\pi:\sqC\to \sigma_P,\]
  where a fiber of $\pi$ is an enhanced dual graph with certain edge lengths in $\R_{\ge 0}$ and vertices weighted by genera of the corresponding components.
  \item If $S$ is a general logarithmic scheme, the above construction can be globalized. Since this is not required in this paper, we omit the details here. Interested readers may consult~\cite{trostack}.
\end{itemize}

\subsection{Relative geometry and expanded degenerations}
We briefly review the technique of expansions of maps to an SNC target pair, introduced by Ranganathan in~\cite{expansion}. Let $X$ be a smooth projective variety and $D\subset X$ a simple normal crossings divisor. Consider a logarithmic stable map $C\to(X,D)$ over $\Spec(\N\to k)$. This induces a map on the tropicalization
$$\sqC\to \Sigma_{(X,D)}.$$
In general, this tropical map does not map cones surjectively to cones, and is therefore not \textit{combinatorially transverse}. To fix it, one can first choose a polyhedral subdivision of the target $\widetilde{\Sigma}\to \Sigma$. Once this is done, the map $\sqC\to \widetilde{\Sigma}$ may fail to be polyhedral, and it is necessary to do a further modification $\widetilde{\sqC}\to\sqC$ to obtain a combinatorial transverse map $\widetilde{\sqC}\to\widetilde{\Sigma}$. Tropically we have a Cartesian diagram:
\begin{center}
  \begin{tikzcd}
    \widetilde{\sqC}\arrow[r]\arrow[d] \arrow[dr, phantom, "\square"]& \widetilde{\Sigma}\arrow[d]\\
    \sqC\arrow[r] & \Sigma.
  \end{tikzcd}
\end{center}
Since both $\widetilde{\Sigma}\to \Sigma$ and $\widetilde{\sqC}\to\sqC$ are polyhedral subdivisions, they correspond to some logarithmic modifications $\widetilde{C}\to C$ and $\widetilde{X}\to X$. And the resulting map $\widetilde{C}\to \widetilde{X}$ would have equidimensional fibers.

The construction is globalized by Ranganathan in~\cite{expansion}, and is a generalization of the smooth pair case proposed by Kim in~\cite{Kim}. The main theorem is the following.
\begin{thm}
  There exists a proper moduli space $K_\Gamma(X)$ of logarithmic stable maps to expansions of $(X,D)$ of fixed contact orders $\Gamma$ with flat source and target families, such that the universal curve transverse to the universal target at every point on the base.
\end{thm}

\subsection{Elliptic singularities}
We review the definition and properties of elliptic singularities, which is a central object of study for \textit{radially aligned} logarithmic curves. A key observation was written by Vakil in~\cite[Lemma 5.9]{vakthesis} and is later on explored by many others in more details, e.g. in~\cite{smyth}.

Let $C$ be a reduced curve, $p\in C$ an isolated singular point, and $\pi:\widetilde{C}\to C$ the normalization of $C$ at $p$. The \textit{genus} of the singular point $p$ is defined as
$$g=\dim(\pi_*(\mathcal{O}_{\widetilde C}/\mathcal{O}_C))-|\pi^{-1}(p)|+1.$$

Roughly speaking, $g$ is the number of extra conditions required for a function to descend from $\widetilde{C}$ to $C$ besides from the obvious topological ones, which only ask for functions to agree at points of $\pi^{-1}(p)$.

A curve is \textit{Gorenstein} if the dualizing sheaf of the curve is invertible. Smyth proved in~\cite{smyth} that, for each integer $m\ge 1$, there is a unique genus 1 Gorenstein curve singularity with $m$ branches. They are called \textit{elliptic singularities} and are classified as follows.

\begin{defn}
  The singular point $p$ is an \textbf{elliptic }$\boldsymbol{m}$\textbf{-fold point} of $C$ if
  $$\hat{\mathcal{O}}_{C,p}=\begin{cases}
    k[[x,y]]/(y^2-x^3) &m=1\text{ (ordinary cusp)},\\
    k[[x,y]]/(y^2-x^2y) &m=2\text{ (ordinary tacnode)},\\
    k[[x_1,\dots,x_{m-1}]]/I_m &m\ge 3\text{ (} m\text{ general lines through the origin in }\A^{m-1}\text{)},\\
  \end{cases}$$
  where $I_m=(x_hx_i-x_hx_j:i,j,h\in\{1,2,\dots,m-1\}$ distinct$)$.
\end{defn}
\begin{ex}
Consider a function on the parabola $y=x^2$ and a function on the line $y=0$. If they are to descend to the tacnode $y^2=x^2y$, one extra condition asking them to have the same derivative at origin is required. Therefore an ordinary tacnode is a genus 1 singularity.
\end{ex}

An elliptic $m$-fold singularity is formed when we contract a genus 1 component with $m$ external nodes in a smoothing family. For example, suppose we have a smoothing family of a nodal curve $C_0$, where $C_0$ has an irreducible genus 1 component $E$ and two irreducible genus 0 components each intersecting $E$ at a unique nodal point. If $E$ is to be contracted, we are forced to replace $E$ by a genus 1 singularity with $m=2$ branches, which must be a tacnode.

\begin{figure}[h]
  \centering
  \captionsetup{justification=centering,margin=2cm}
    \begin{tikzpicture}[scale=1]
      \draw (0,0) ellipse (1.5cm and 1cm);
      \draw (-0.6,0.08) arc (-130:-50:1);
      \draw (0.35,-0.06) arc (50:130:0.5);
      \draw (2,1.23) circle (1);
      \draw (2,-1.23) circle (1);
      \draw[dashed] (3,1.23) arc (0:180:1 and 0.3);
      \draw (3,1.23) arc (0:-180:1 and 0.3);
      \draw [dashed](3,-1.23) arc (0:180:1 and 0.3);
      \draw (3,-1.23) arc (0:-180:1 and 0.3);
      \filldraw (1.26,0.54) circle (1pt);
      \filldraw (1.26,-0.54) circle (1pt);
      \node at (-1.7,1) {$E$};
      \node at (3.6,1.5) {$\P^1$};
      \node at (3.6,-1.5) {$\P^1$};
      \node at (1,-3) {$C_0$};
      \node at (8,-3) {tacnode};
      \draw [->,
        line join=round,
        decorate, decoration={
            zigzag,
            segment length=4,
            amplitude=0.4,post=lineto,
            post length=2pt
        }]  (4,0) -- (6,0);
      \draw [smooth,samples=100,domain=7:9] plot(\x,{((\x-8)*(\x-8)*1-0.3)});
      \draw (7,-0.3) -- (9,-0.3);
      \end{tikzpicture}
    \caption{A tacnode singularity in a contracting family.}
    \label{ellipticsing}
\end{figure}
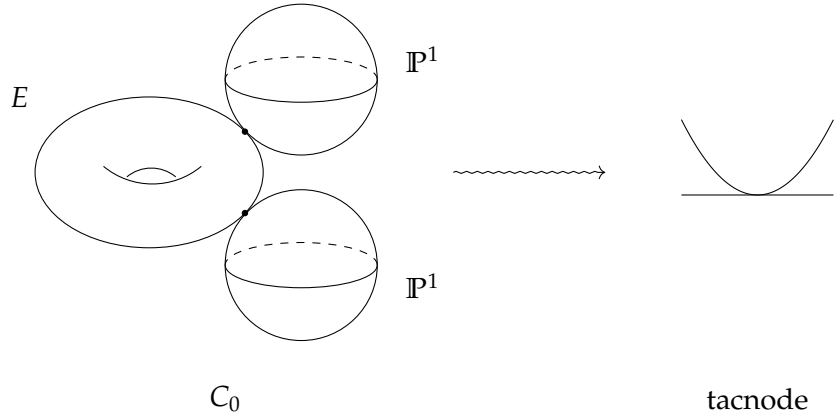

\subsection{Radially aligned curves}\label{sec:radcuve}
We consider logarithmic curves of genus 1 that are \textit{radially aligned}. The data in a radial aligned curve suffices to determine the type of elliptic singularities if the elliptic component is to be contracted. We make this precise. We refer readers to~\cite{RSW1} for further details and proofs stated in this section.

Let $C\to S$ be a logarithmic curve of genus 1 with tropicalization $\sqC$. The \textit{circuit} of $C$ is the minimal genus 1 subcurve, and analogously the circuit of $\sqC$ is the union of vertices whose complement contains no component of genus 1.

There is a piecewise linear function $\lambda$ of distance from each vertex $v$ to the circuit of $\sqC$. This is done by summing over all the edge lengths $$\lambda(v)=\sum_{i=1}^k l(e_i)\in \overline{M}_S$$ of edges in the unique path from $v$ to the circuit. The path is unique because the tropical curve is essentially a tree if we treat the circuit as a single point of destination.

\begin{defn}
  A genus 1 logarithmic curve $C\to S$ is \textbf{radially aligned} if for all geometric point $s\in S$ and vertices $v,w$ of $\sqC_s$, the function $\lambda$ satisfies that $\lambda(v)$ and $\lambda(w)$ are always comparable.
\end{defn}

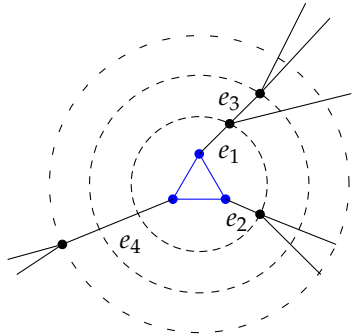
\begin{figure}[h]
  \centering
  \captionsetup{justification=raggedright}
  \begin{tikzpicture}[scale=1]
    \tikzstyle{every node}=[font=\small]
    \filldraw [blue](0,1*0.4) circle (1.6pt);
    \filldraw [blue](0.866*0.4,-0.5*0.4) circle (1.6pt);
    \filldraw [blue](-0.866*0.4,-0.5*0.4) circle (1.6pt);
    \draw [blue](0,1*0.4)--(0.866*0.4,-0.5*0.4)--(-0.866*0.4,-0.5*0.4)--(0,1*0.4);
    \filldraw(1*0.4,2*0.4) circle (1.6pt);
    \filldraw(2*0.4,-1*0.4) circle (1.6pt);
    \filldraw(2*0.4,3*0.4) circle (1.6pt);
    \filldraw(-4.5*0.4,-2*0.4) circle (1.6pt);
    \draw (0,1*0.4)--(1*0.4,2*0.4)--(2*0.4,3*0.4)--(3.5*0.4,6*0.4);
    \draw (2*0.4,3*0.4)--(4.3*0.4,5.5*0.4);
    \draw (1*0.4,2*0.4)--(5*0.4,3*0.4);
    \draw (-0.866*0.4,-0.5*0.4)--(-4.5*0.4,-2*0.4)--(-6.3*0.4,-2.5*0.4);
    \draw (-4.5*0.4,-2*0.4)--(-6*0.4,-3*0.4);
    \draw (0.866*0.4,-0.5*0.4)--(2*0.4,-1*0.4)--(4.5*0.4,-2*0.4);
    \draw (2*0.4,-1*0.4)--(4*0.4,-3*0.4);
    \draw [densely dashed] (0,0) circle (2.236*0.4);
    \draw [dashed] (0,0) circle (3.606*0.4);
    \draw [loosely dashed] (0,0) circle (4.924*0.4);
    \node at (0.4,0.4) {$e_1$};
    \node at (0.5,-0.5) {$e_2$};
    \node at (0.4,1.1) {$e_3$};
    \node at (-0.9,-0.8) {$e_4$};
  \end{tikzpicture}
  \caption{Tropicalization of a genus 1 logarithmic curve. The core is the triangle in the middle. The curve is radially aligned if $e_1=e_2$ and $e_1+e_3<e_4$, in which case the distances from each vertex to the core are totally ordered.}
  \label{radcurve}
\end{figure}

Equivalently, radially alignment asks for a total ordering on the vertices that are outside the core of $\sqC$. See Figure \ref{radcurve} for an example.

This example also suggests that the radial alignment condition locally is a toric modification arising from the subdivision of the tropical moduli by $e_1=e_2$ and $e_1+e_3=e_4$. This is called a \textit{logarithmic modification}.

\begin{defn}
  A morphism of logarithmic schemes $X'\to X$ is a \textbf{logarithmic modification} if it is locally pulled back from a toric modification.
\end{defn}

Let $\mathfrak M_{1,n}$ be the logarithmic stack of $n$-pointed, genus 1 pre-stable curves. Let $\mathfrak M^{\rad}_{1,n}$ be the category fibered in groupoids over logarithmic schemes whose fiber over $S$ is the groupoid of radially aligned  logarithmic curves of genus 1 and $n$ marked points. We have the following result.
\begin{thm}
  $\mathfrak M^{\rad}_{1,n}$ is a logarithmic modification of $\mathfrak M_{1,n}$.
\end{thm}

One consequence of the theorem is that the stack $\mathfrak M^{\rad}_{1,n}$ is logarithmically smooth, because it is logarithmically \etale over the logarithmically smooth stack $\mathfrak M_{1,n}$. However, Ranganathan--Santos-Parker--Wise proved further that the logarithmic structure of $\mathfrak M^{\rad}_{1,n}$ is locally free, meaning it is in fact a smooth algebraic stack~\cite[Corollary~3.4.3]{RSW1}.

Finally, to produce an elliptic singularity, we need to pick a distance and consider the circle from the circuit of the tropicalization of radius equal to that distance, and perform a contraction of everything within the circle. We explain this process with the example in Figure \ref{radcurve} again.

\begin{itemize}
  \item Pick an integer $m$ smaller than the number of marked points $n=7$, e.g. $m=5$.
  \item There is a smallest radius $\delta$ such that the circle of radius $\delta$ from the core has inner valence $\le m=5$ and outer valence $\ge m=5$. In this example, $\delta=e_1+e_3$.
  \item Introduce vertices if possible on the intersections of the circle and the edges. This is a partial destabilization of the curve $\widetilde{C}\to C$.
  \item Perform the contraction of the circle of radius $\delta$ to get $\widetilde{C}\to \overline{C}$, and $\overline{C}$ is a curve with a no-worse-than $m$-fold elliptic singularity.
\end{itemize}

The precise statement can be found at~\cite[Theorem 3.7.1]{RSW1}, which involves an extra notion of $m$-stability.

\subsection{The factorization condition on stable maps}\label{sec:factorization}
Let $\overline{\mathcal M}_{1,n}(Y,\beta)$ be the moduli space of $n$-pointed, genus 1 logarithmic stable maps to $Y$ with image class $\beta$. There is a forgetful morphism by forgetting the map structure
$$\overline{\mathcal M}_{1,n}(Y,\beta)\to \mathfrak M_{1,n}.$$

We define the moduli space $\widetilde{\VZ}_{1,n}(Y,\beta)$ of radially aligned maps by the following Cartesian diagram:

\begin{center}
  \begin{tikzcd}
    \widetilde{\VZ}_{1,n}(Y,\beta)\arrow[r]\arrow[d]\arrow[dr, phantom, "\square"] &\overline{\mathcal M}_{1,n}(Y,\beta)\arrow[d]\\
    \mathfrak{M}_{1,n}^{\text{rad}} \arrow[r] & \mathfrak{M}_{1,n}.
  \end{tikzcd}
\end{center}
The space $\widetilde{\VZ}_{1,n}(Y,\beta)$ parametrizes logarithmic stable maps from radial genus 1 and $n$-pointed logarithmic curves to $Y$ with class $\beta$. The map $$\widetilde{\VZ}_{1,n}(Y,\beta)\to\overline{\mathcal M}_{1,n}(Y,\beta)$$ is a logarithmic modification because it is the pullback of a logarithmic modification.

Let $C/S$ be a family of radial curves and $\lambda$ be the piecewise-linear function of distance to the circuit. Let $C\to Y$ be a stable map in $\widetilde{\VZ}_{1,n}(Y,\beta)$. There is a section $\delta$ parameterizing the \textit{contraction radius}, where over $s\in S$, the radius $\delta_s$ is the smallest $\lambda(v)$ for a non-contracted component $v\in C_s$. This section $\delta$, similar to the case in the concluding example in Section \ref{sec:radcuve}, induces a partial destabilization $\widetilde{C}\to C$ and a contraction $\widetilde{C}\to \overline{C}$.

\begin{defn}\label{def:faccon}
  A stable map $[f:C\to Y]\in\widetilde{\VZ}_{1,n}(Y,\beta)$ satisfies the \textbf{factorization property} if the composition $\widetilde{C}\to C\to Y$ factors through $\overline{C}$:
  \begin{center}
    \begin{tikzcd}
      \widetilde{C}\arrow[r]\arrow[dr] & C \arrow[r] & Y\\
      & \overline{C}\arrow[ur].&
    \end{tikzcd}
  \end{center}
\end{defn}

Write ${\VZ}_{1,n}(Y,\beta)$ for the substack of maps satisfying the factorization property. We list two important properties for these spaces, which are proved in~\cite{RSW1}.

\begin{thm}
  If $Y$ is proper, then ${\VZ}_{1,n}(Y,\beta)$ is proper.
\end{thm}

\begin{thm}\label{thm:abssmooth}
  The stack ${\VZ}_{1,n}(\P^r,d)$ is smooth and is of expected dimension.
\end{thm}

Theorem \ref{thm:abssmooth} is proved by showing that the map to the stack of radial curves is unobstructed. In particular, because the resulting space is smooth, this construction provides another point of view of the Vakil--Zinger desingularization for the space of genus 1 stable maps to projective spaces.

\section{Construction of moduli for a subtoric divisor of $\P^{\underline n}$}\label{sec:moduli}

\subsection{Notation and setup}

Let $\P^{\underline{n}}=\P^{n_1}\times\cdots\times\P^{n_a}$ be a product of projective spaces and $D\subset \partial \P^{\underline{n}}$ a (possibly empty) subset of the toric boundary divisors. We construct a moduli space of genus 1 radial stable maps to $\P^{\underline n}$ relative to $D$ that are \textit{well-spaced}. We also investigate the behaviors of the tropicalization of those stable maps. We will prove in Section \ref{sec:etale} that the moduli space constructed is logarithmically smooth.

We start with the space $\widetilde{\VZ}{}^\Gamma_{1,n}(\P^{\underline{n}},D)$ of logarithmic stable maps from radially aligned $n$-pointed genus 1 curves to expansions of $(\P^{\underline{n}},D)$ with contact orders $\Gamma$.

Let $f:C\to \P^{\underline{n}}$ be such a map to $\P^{\underline{n}}$. If the genus 1 component of $C$ is not contracted, then it is straightforward to show that the map is unobstructed, and thus the corresponding stratum in the moduli space will have expected dimension.

\subsection{Well-spaced maps to $(\P^{\underline n},\emptyset)$}\label{sec:desing}
We begin by defining the well-spacedness condition for maps to a product of projective spaces relative to the empty divisor. Let $f:C\to \P^{\underline n}$ be a map from a radial genus 1 curve. It comes with canonical projections
$$f_i:C\to\P^{\underline n}\to \P^{n_i}.$$
In particular, there is degree information on components of C for each $f_i$ that determines whether a component is contracted or not. Similar to the process described in Section \ref{sec:factorization}, we could ask all $f_i$ to satisfy the factorization condition in Definition \ref{def:faccon}.

\begin{defn}\label{def:well-empty}
  A map $[f:C\to \P^{\underline n}]\in\widetilde{\VZ}_{1,n}(\P^{\underline n},\boldsymbol{d})$ satisfies the \textbf{well-spacedness condition} if $f$ and all projections $f_i:C\to \P^{\underline n}\to \P^{n_i}$ satisfy the factorization condition.
\end{defn}
Write $\VZ_{1,n}(\P^{\underline n},\emptyset)$ for the substack of maps satisfying the well-spacedness condition.
\begin{ex}
  Consider a map $f:C\to \P^{n_1}\times \P^{n_2}$ where the tropicalization of $f$ is the one drawn in Figure \ref{fig:example1}, with decoration of the bidegree information on each vertex for the two projections $f_i:C\to \P^{n_i}$.

  \begin{figure}[h]
    \centering
    \captionsetup{justification=raggedright}
    \begin{tikzpicture}[scale=1]
      \tikzstyle{every node}=[font=\small]
      \filldraw [blue](0,0) circle (1.6pt);
      \node at (0,0.3) {$g=1,(0,0)$};
      \filldraw(-0.8,-0.6) circle (1.6pt);
      \node at (-1.4,-0.6) {$(0,0)$};
      \filldraw(0.8,-1) circle (1.6pt);
      \node at (0.2,-1) {$(1,0)$};
      \filldraw(0.6,-2) circle (1.6pt);
      \node at (0.2,-1.7) {$(1,0)$};
      \filldraw(1.5,-1.5) circle (1.6pt);
      \node at (2,-1.75) {$(1,0)$};
      \filldraw(-2,-2.3) circle (1.6pt);
      \node at (-2.7,-2.3) {$(0,1)$};
      \filldraw(-1,-2.3) circle (1.6pt);
      \node at (-0.5,-2.1) {$(0,1)$};
      \draw (0,0)--(-0.8,-0.6)--(-2,-2.3)--(-1.7,-3);
      \draw (-2,-2.3)--(-2.3,-3);
      \draw (-0.8,-0.6)--(-1,-2.3)--(-0.7,-3);
      \draw (-1,-2.3)--(-1.3,-3);
      \draw (0,0)--(0.8,-1)--(0.6,-2)--(0.3,-3);
      \draw (0.6,-2)--(0.9,-3);
      \draw (0.8,-1)--(1.5,-1.5)--(1.2,-3);
      \draw (1.5,-1.5)--(1.8,-3);
      \draw [->] (3,0)--(3,-3);
      \node at (4,0.3) {distance from the circuit};
      \filldraw(3,0) circle (1pt);
      \filldraw(3,-0.6) circle (1pt);
      \filldraw(3,-1) circle (1pt);
      \filldraw(3,-1.5) circle (1pt);
      \filldraw(3,-2) circle (1pt);
      \filldraw(3,-2.3) circle (1pt);
      \draw [loosely dashed] (-0.8,-0.6)--(3,-0.6);
      \node at (3.3,-0.6) {$e_1$};
      \draw [loosely dashed] (0.8,-1)--(3,-1);
      \node at (3.3,-1) {$e_2$};
      \draw [loosely dashed] (1.5,-1.5)--(3,-1.5);
      \node at (3.65,-1.5) {$e_2+e_3$};
      \draw [loosely dashed] (0.6,-2)--(3,-2);
      \node at (3.65,-2) {$e_2+e_4$};
      \draw [loosely dashed] (-2,-2.3)--(3,-2.3);
      \node at (3.65,-2.4) {$e_1+e_5$};
    \end{tikzpicture}
    \caption{A decorated genus 1 radially aligned tropical curve.}
    \label{fig:example1}
  \end{figure}
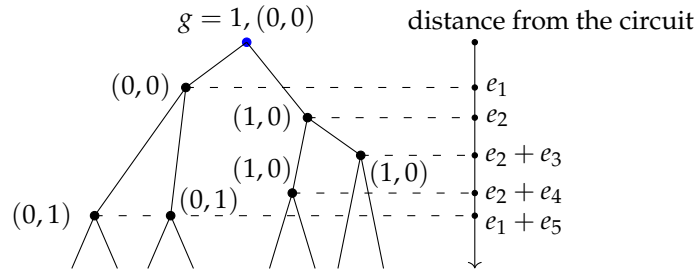

  The factorization property for $f:C\to \P^{n_1}\times \P^{n_2}$ requires $f$ to factor through the contraction of the curve of radius $e_1$, for $f_1:C\to \P^{n_1}$ the contraction of radius also $e_1$, and for $f_2:C\to \P^{n_2}$ the contraction of radius $e_1+e_5$. The map $f$ is well-spaced if the factorization property holds for all three maps.
\end{ex}

As shown by the above example, the well-spacedness condition in Definition \ref{def:well-empty} is an explicit requirement that will (eventually) desingularizes for the space ${\VZ}_{1,n}(\P^{\underline n},\boldsymbol{d})$.

\subsection{The logarithmic multiplicative group $\Glog$}
We will make use of the logarithmic multiplicative group $\Glog$ to define the well-spacedness condition when $D\not=\emptyset$, which we develop in this section. For more details on $\Glog$, see~\cite{Molcho} or~\cite{dhruv-wise}.

\begin{defn}
  The logarithmic multiplicative group $\Glog$ is the functor $\textbf{LogSch}^\text{op}\to\textbf{Sets}$ defined by
  $$\Glog(X):=\Gamma(X,\Mgp_X).$$

\end{defn}
For example, if $X$ has a trivial logarithmic structure, then $\Glog(X)=\G_m(S)$. So morally $\Glog$ partially compactifies $\G_m$. However, $\Glog$ is not representable in the category of logarithmic schemes, as shown in~\cite[Lemma 2.2.7.3]{Molcho}.

\begin{prop}\label{prop:glog}
  The map $\P^1\to\Glog$ given locally by $(x,y)\mapsto x^{-1}y$ is a logarithmic modification, where $\P^1$ is with its toric logarithmic structure.
\end{prop}

Since logarithmic modifications are logarithmically \etale, $\Glog$ has a logarithmically \etale cover by a proper logarithmic scheme $\P^1$. Therefore $\Glog$ should be seen as proper. However, the only schematic compactification of $\G_m$ is $\P^1$, which is not a group object in the category of schemes. So $\Glog$ is not representable in the category of logarithmic schemes.

The following result is a direct generalization of Prosition \ref{prop:glog}.

\begin{prop}\label{prop:glogtoric}
  $\Glog^r$ admits a logarithmically \etale cover by any complete toric variety of dimension~$r$.
\end{prop}

\subsection{Expansion and $\Glog^r$-bundle} Suppose now $D\not=\emptyset$, and the circuit of $C$ is contracted under the map $f$. If the circuit is contracted to the interior $\P^{\underline n}\setminus D$, then locally around the circuit, the map $f$ behaves like a map to $(\P^{\underline n},\emptyset)$, and we could impose the well-spacedness condition defined in Definition \ref{def:well-empty}.

On the other hand, if the circuit is contracted into $D$, it is necessary to consider the map to the expanded target. Write $W$ for stratum of $D$ that is of the highest codimension and which the circuit is contracted into, and set $r$ to be the codimension of $W$ in $\P^{\underline n}$. Notice that $W$ itself is some product of projective spaces. We demonstrate that we could treat $f:C\to \P^{\underline{n}}$ locally as a map to a $\Glog^r$-bundle over $W$.

Consider first that $D\subset \P^{\underline n}$ is irreducible and smooth. Then in the setting above, the genus 1 component is contracted into $W=D$ with $r=1$. The expansion of $\P^{\underline n}$ along $W$ adds potentially multiple copies of the projective completion of the normal bundle $\P(\mathcal{N}_W\oplus 1)$ to $\P^{\underline n}$. In particular, the highest level is a complete toric variety. Similarly, if $W$ has codimension $r$, then the expanded degeneration of $\P^{\underline n}$ gives a $(\P^1)^r$-bundle over $W$.

Now we apply Proposition \ref{prop:glogtoric}, which says that any complete toric variety of dimension~$r$ admits a map to $\Glog^r$.

\begin{lem}\label{lem:maptoglog}
  Let $C'\subset C$ be the maximal genus 1 subcurve contracted into $W$ of codimension~$r$. Then locally there is a map $C'\to B$, where $B$ is a $\Glog^r$-bundle over $W$.\qed
\end{lem}

\begin{rem}
  An equivalent approach is to begin with the $(\P^1)^r$-bundle and remove all the 0 and $\infty$ sections to get a $\G_m^r$-bundle. Since $\Glog$ compactifies $\G_m$, we could compactify this map to get $C'\to B$ for a $\Glog^r$-bundle $B/W$.
\end{rem}

\begin{ex}
  Lemma \ref{lem:maptoglog} is very visible in the case where the target is the smooth pair $(\P^1,D=\text{pt})$. The projective completion of $D$ is $\P(\mathcal{N}_D\oplus 1)=\P^1$, and the $k$-fold expansion is $\P^1[k]=\P^1\cup\dots\cup \P^1$, where $C'$ will map to the highest level $\P^1$. Also $\P^1$ admits a map to $\Glog$ by Proposition 3.2, and this is the map we want. Alternatively, we could remove 0 and $\infty$ to get a map to $\G_m$, and compactify again to obtain a map to $\Glog$.
\end{ex}

\subsection{The well-spacedness condition}
Once we have the map $C'$ to $B/W$, which is a $\Glog^r$-bundle over $W$ from Lemma \ref{lem:maptoglog}, we can use it to define the well-spacedness condition.
\begin{defn}\label{def:ws}
  A map $f\in \widetilde{\VZ}{}^\Gamma_{1,n}(\P^{\underline{n}},D)$ is \textbf{well-spaced} if:
  \begin{itemize}
    \item either the genus 1 component is not contracted, or
    \item all maps $C'\to B\to QB$, where $QB$ is any quotient bundle, satisfy the factorization property, and the map to the base $C'\to W$ satisfy the well-spacedness condition defined in Definition \ref{def:well-empty}.
  \end{itemize}
\end{defn}
This generalizes Definition \ref{def:well-empty} because if $D=\emptyset$, then the bundle $B=W\times \Glog^0$ is trivial with $r=0$, and we are simply asking the map to $W=\P^{\underline n}$ and all projections to satisfy the factorization condition.

Let ${\VZ}^\Gamma_{1,n}(\P^{\underline{n}},D)$ be the substack of $\widetilde{\VZ}{}^\Gamma_{1,n}(\P^{\underline{n}},D)$ of well-spaced maps.

\begin{prop}
 ${\VZ}^\Gamma_{1,n}(\P^{\underline{n}},D)$ is proper.
\end{prop}
This is immediate because of~\cite[Proposition 3.4.5]{RSW2}, which says that the factorization property satisfies the valuative criterion for properness.

This definition of well-spacedness generalizes the notion of \textit{well-spaced maps to toric varieties} defined in~\cite{RSW2}. For a map $f:C\to Z$ to a toric variety $Z$ with character lattice $N^\vee$, Ranganathan requires the map $$N^{\vee}_{T/H}\to N^\vee\to \Gamma(C,M_C^{\text{gp}})$$
to satisfy the factorization property for any subtorus $H$ of the dense torus $T\subset Z$. This is equivalent to asking $$C\to N\otimes \Glog\to N_{T/H}\otimes \Glog$$ to factor through some $\overline{C}_H$. Our quotient bundle $QB$ plays the role of $N_{T/H}\otimes \Glog$ here.

\subsection{Tropical well-spacedness condition and moduli}\label{sec:tropmoduli}

We define the \textit{tropical well-spacedness condition} and construct the moduli space of tropical maps that are well-spaced. This moduli space is a generalized cone complex, with each cone parameterizing a certain stratum in $\VZ^\Gamma_{1,n}(\P^{\underline n},D)$.

First we recall the definition of tropical well-spacedness condition for maps to a toric fan, defined by Ranagathan in~\cite{RSW2}.

\begin{defn}
  A tropical genus 1 stable map $\sqC\to \Sigma$ to a toric fan $\Sigma$ is \textbf{tropically well-spaced for toric targets} if for each character $\chi:\N_\R\to \R$, the induced map
  $$F:\sqC\to \R$$
  satisfies:
  \begin{itemize}
    \item either no open neighborhood of the circuit $\sqC_0$ of $\sqC$ is contracted, or
    \item let $t_1,\dots,t_k$ be the flags whose base is mapped to $F(\sqC_0)$ but along which $F$ has non-zero slope; then the minimum of the distances $\{d(t_i,F(\sqC_0))\}_{i=1}^k$ occurs at least three times.
  \end{itemize}
\end{defn}

In our case, the use of $\Glog^r$ helps the situation because it almost has a ``fan'' which looks like an entire copy of $\R^r$. We could then analyze whether tropical maps to $\Sigma_{(\P^{\underline n},D)}$ are well-spaced by the following procedure.

Fix a radial genus 1 stable map $f:C\to \P^{\underline n}$. Its tropicalization is a map
$$w:\sqC\to \Sigma_{(\P^{\underline n},D)}$$
of cone complexes from a genus 1 tropical curve to the tropicalization of $(\P^{\underline n},D)$. We assume $w$ is combinatorially transverse. Note that for $w$, each vertex of $\sqC$ has the information of a \textit{degree} associated to it, corresponding to the degree of the component under $f$.

\begin{itemize}
  \item If the circuit of $\sqC$ has non-zero degree, then the map $C'\to B$ to the $\Glog^r$-bundle does not contract the circuit. Therefore this map $f$ is well-spaced.
  \item If the circuit has degree 0, then by Definition \ref{def:ws}, we need the map $C'\to B$ to satisfy the factorization condition for any quotient bundle. This is similar to asking a map to toric target to be well-spaced, and we use the tropical well-spacedness condition defined above.
\end{itemize}

\begin{defn}
  A combinatorially transverse genus 1 tropical map $w:\sqC\to \widetilde{\Sigma}_{(\P^{\underline n},D)}$ from a radial genus 1 curve is \textbf{tropically well-spaced} if
  \begin{itemize}
    \item either the circuit has non-zero degree, or
    \item the induced map $\sqC'\to \text{Trop}(B)=\text{Trop}(W)\times \R^r\to \R^r$ satisfies the \textit{tropical well-spacedness condition for toric targets}.
  \end{itemize}
\end{defn}

Denote by $W^\Gamma_{1,n}(\P^{\underline n},D)$ the moduli space of well-spaced tropical maps. Below is a brief summary of how to explicitly construct this tropical moduli space.

\begin{itemize}[label={},leftmargin=*]
  \item \textbf{Step 1.} We start with the tropical moduli space of genus one stable maps to $(\P^{\underline n}, D)$. A modification to source curves is required to make them radial aligned (see Section \ref{sec:radcuve}).\\

  \item \textbf{Step 2.} For each tropical stable map $\sqC\to \Sigma_{(\P^{\underline n},D)}$, perform a polyhedral subdivision on the target and pull back the construction to the source. The resulting tropical map $\widetilde{\sqC}\to\widetilde{\Sigma}_{(\P^{\underline n},D)}$ will be transverse by~\cite{expansion}. At this stage, the resulting tropical moduli space will be the tropicalization of $\widetilde{\VZ}{}^\Gamma_{1,n}(\P^{\underline{n}},D)$, which we write as $\widetilde{W}{}^\Gamma_{1,n}(\P^{\underline{n}},D)$.\\

  \item \textbf{Step 3.} Take the subcomplex $W^\Gamma_{1,n}(\P^{\underline n},D)$ of $\widetilde{W}{}^\Gamma_{1,n}(\P^{\underline{n}},D)$ of tropical well-spaced maps.
\end{itemize}

\section{Logarithmic smoothness of the moduli spaces}\label{sec:etale}
This section analyzes the behavior of the moduli spaces $\VZ^\Gamma_{1,n}(\P^{\underline n},D)$ when we add an additional divisor $D_0$, subject to some generality conditions, to the divisor $D$ of the target space $\P^{\underline n}$. The map forgetting the divisor $D_0$ locally around a generic lift is \etale (see Proposition \ref{mainlem} for the precise statement), which we will later refer to as Property ($\star$) for convenience. By sequentially adding in all divisors $\partial \P^{\underline n}\setminus D$ not in the toric boundary, we could compare stable maps to $(\P^{\underline n},D)$ with maps to the toric $\P^{\underline n}$. We use this technique to prove that the moduli spaces $\VZ^\Gamma_{1,n}(\P^{\underline n},D)$ are logarithmically smooth.
\subsection{Adding one coordinate divisor}
We call the pullback of any hyperplane under any projection $\P^{\underline n}=\P^{n_1}\times\cdots\times\P^{n_a}\to \P^{n_i}$ a \textit{coordinate divisor}. The toric boundary divisors of $\P^{\underline n}$ consist of pullbacks of all coordinate hyperplanes of $\P^{n_i}$ under the projections, and $D$ is a subset of the toric boundary divisors.

A generic coordinate divisor $D_0$ different from all divisors in $D$ would intersect a given curve transversely in a number of points. These new contact orders can be used to build a new moduli space. Proposition \ref{mainlem} is the central result of the section.

\begin{prop}\label{mainlem}
Given a generically chosen coordinate divisor $D_0$ of $\P^{\underline{n}}$ and a map $f\in \VZ^{\Gamma}_{1,n}(\P^{\underline{n}},D)$, there exists a lift $\widetilde{f}\in \VZ^{\Gamma'}_{1,n+m}(\P^{\underline{n}},D+D_0)$ where
\begin{enumerate}
  \item the contact orders for the original $n$ marked points with $D_0$ are 0;
  \item the contact orders for the new $m$ points with $D_0$ are 1, and with any other components of $D$ are 0.
\end{enumerate}
Also the forgetful map
$$\VZ^{\Gamma'}_{1,n+m}(\P^{\underline{n}},D+D_0)\to \VZ^{\Gamma}_{1,n}(\P^{\underline{n}},D),$$
locally around $\widetilde{f}$ and $f$ is \etale, i.e. there is an isomorphism
$$\hat{\mathcal O}_{\VZ^{\Gamma}_{1,n}(\P^{\underline{n}},D),f}\xrightarrow{\cong} \hat{\mathcal O}_{\VZ^{\Gamma'}_{1,n+m}(\P^{\underline{n}},D+D_0),\widetilde{f}}$$
between the two complete local rings.
\end{prop}

\begin{note}
  The number $m$ is uniquely determined as the sum of degrees of those components that a generic $D_0$ will intersect. See Example \ref{ex:genericonce} for an illustration.
\end{note}

We demonstrate the existence of $\widetilde{f}$ first. This is immediate because the new divisor $D_0$ is generic.

\begin{lem}\label{transint}
  The set of coordinate divisors of $\P^{\underline n}$ which intersect a fixed curve $C$ transversely at points outside $\{p_1,\dots,p_n\}$ is open among all coordinate divisors.
\end{lem}
\begin{pf}
  In terms of the coefficients for the defining equation of the coordinate divisor, both intersecting $C$ with higher tangency and intersecting $C$ at the given points $p_i$ are closed conditions. Therefore the complement that we require is open.\qed
\end{pf}

Thus the same underlying map $f$ can be seen as an element $\widetilde{f}\in \VZ^{\Gamma'}_{1,n+m}(\P^{\underline{n}},D+D_0)$, subject to a choice of labeling of the new $m$ marked points. Note that if $f$ is well-spaced, so will $\widetilde{f}$.

To show the forgetful map is \etale, we use the infinitesimal lifting criterion. Fix an infinitesimal logarithmic deformation of $f$ over $k[\epsilon]/(\epsilon^2)$. We show the same underlying family can also be seen as a deformation in $\VZ^{\Gamma'}_{1,n+m}(\P^{\underline{n}},D+D_0)$. For this to happen, the family must preserve the contact orders with both $D$ and $D_0$. The contact orders with $D$ is preserved by being a deformation in $\VZ^{\Gamma}_{1,n}(\P^{\underline{n}},D)$, so it suffices to check the same for $D_0$.

\begin{lem}\label{famint}
  Let $F$ be an infinitesimal logarithmic deformation of $f:C\to \P^{\underline n}$ with $F_0=f$, and $D_0$ a generically chosen coordinate divisor as above. Then $F_\epsilon$ intersects $D_0$ with contact orders $(1,1,\dots,1)$.
\end{lem}

\begin{pf}
  We check the intersection at source level. The equation $h$ for $D_0$ is a section of the line bundle $\O(1)$ on $\P^{n_1}$ after a change of coordinates. The line bundle $\O(1)$ pulls back to a line bundle $\mathcal{L}=f^*\O(1)$ over the source $C$, and $h$ pulls back to a section $f^*h$:
  \begin{center}
    \begin{tikzcd}
      \mathcal{L}=f^*\O(1)\arrow[r]\arrow[d]& \O(1)\arrow[d]\\
      C\arrow[r,"p_1\circ f"]&\P^{n_1}.
    \end{tikzcd}
  \end{center}
  Then $F$ can be realized as an infinitesimal deformation of the 0-section when intersecting $f^*h$ of the line bundle $\mathcal{L}$. Since two sections intersecting transversely is an open condition, the claim of the lemma follows.\qed
\end{pf}
\begin{rem}
  The intuition is that transversality is an open condition both for the coordinate divisor and for the curve, and therefore perturbing the curve should not change the transverse contact orders. The result does not hold if $F_0$ has an intersection point with $D_0$ of order 2, as order 2 is a closed condition and can change to (1,1) when going from the $F_0$ to the generic fiber $F_\epsilon$.
\end{rem}

The above discussion essentially proves Proposition \ref{mainlem}.
\begin{pf}[Proposition \ref{mainlem}]
  By Lemma \ref{transint}, the map $f$ can be lifted to a map $\widetilde{f} \in \VZ^{\Gamma'}_{1,n+m}(\P^{\underline{n}},D+D_0)$ given a choice of the labeling of the $m$ points. Fix one of the labelings, and we show that the map locally around $\widetilde{f}$ and $f$ is \etale by the infinitesimal lifting criterion.

  Take a deformation of $f$ corresponding to a strict square-zero extension. The underlying family is also a deformation of $\widetilde{f} \in \VZ^{\Gamma'}_{1,n+m}(\P^{\underline{n}},D+D_0)$. So a lift of the extension exists. Because the labeling for the $m$ points is fixed, this lift is also unique. Therefore the map $\VZ^{\Gamma'}_{1,n+m}(\P^{\underline{n}},D+D_0) \to \VZ^{\Gamma}_{1,n}(\P^{\underline{n}},D)$ satisfies Property ($\star$).\qed
\end{pf}

\begin{ex}\label{ex:genericonce}
  Consider a map $f$ from a smooth genus 1 curve to $(\P^2,D_0)$ of degree 3, where $D_0=\{x_0=0\}$ is a smooth boundary divisor, with three marked points of contact order 1. Tropically, this maps a vertex with 3 unbounded rays into $\R_{\ge0}$ with the vertex mapping to origin, and 3 rays along the $\R_{>0}$ direction.

  Let $D_1$ be another boundary divisor $D_1=\{x_1=0\}$. Then by Proposition \ref{mainlem}, tropically the corresponding lift will have 3 new rays mapping into $y_{>0}$-axis of the tropicalization $\R_{\ge0}^2$ of $(\P^2,D_0+D_1)$.
  \begin{figure}[h]
    \centering
    \begin{tikzpicture}
      \tikzstyle{every node}=[font=\small]
      \draw [->] (-2,0)--(0,0);
      \draw [blue,->] (-2,0)--(-0.5,0.1);
      \draw [blue,->] (-2,0)--(-0.5,0);
      \draw [blue,->] (-2,0)--(-0.5,-0.1);
      \draw [->] (1,0)--(2.5,0);
      \draw [->] (5,0)--(7,0);
      \draw [->] (5,0)--(5,2);
      \filldraw[blue] (-2,0) circle (1.5pt);
      \filldraw[blue] (5,0) circle (1.5pt);
      \draw [blue,->] (5,0)--(6.5,0.1);
      \draw [blue,->] (5,0)--(6.5,0);
      \draw [blue,->] (5,0)--(6.5,-0.1);
      \draw [blue,->] (5,0)--(5,1.5);
      \draw [blue,->] (5,0)--(5.1,1.5);
      \draw [blue,->] (5,0)--(4.9,1.5);
      \node at (-2,-0.4) [text=blue] {$g=1,\text{degree}=3$};
      \node at (5,-0.4) [text=blue] {$g=1,\text{degree}=3$};
    \end{tikzpicture}
    \caption{Lifting the tropical map when added a new divisor.}
    \label{once}
  \end{figure}
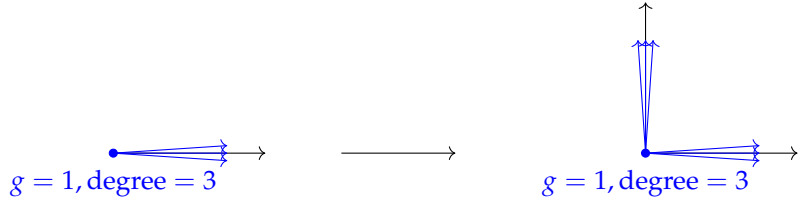
\end{ex}

\subsection{Completing to toric boundary}\label{sec:fulltoric}
We can, by induction, put in extra coordinate divisors to $D$ until it reaches the toric boundary $\partial \P^{\underline{n}}$.

Note that $\P^{\underline n}$ has the property that we can turn any coordinate divisor into a boundary divisor by an automorphism.
\begin{lem}
  Given $D\subset \partial \P^{\underline n}$ and a generic coordinate divisor $D_0$, there exists an automorphism of $\P^{\underline n}$ that fixes $D$ and maps $D_0$ to a toric boundary divisor.
\end{lem}
\begin{pf}
  Without loss of generality, assume $D_0$ is pulled back from a hyperplane $H$ in $\P^{n_1}$. Then pick an automorphism of $\P^{n_1}$ that fixes all the relevant components of $D$ and that maps $H$ to a coordinate hyperplane. This automorphism suffices.\qed
\end{pf}
\begin{rem}
  This may be one of the reasons why this technique can be applied on a product of projective spaces. It might be more difficult to find such an automorphism for a more general target space.
\end{rem}

\begin{cor}\label{cor:fulltoric}
  Given a map $f\in \VZ^\Gamma_{1,n}(\P^{\underline{n}},D)$ and enough generically chosen coordinate divisors, it is possible to lift $f$ to a map $\widetilde{f}\in \VZ^{\Gamma'}_{1,n+m}(\P^{\underline{n}},\partial \P^{\underline{n}})$. Under the forgetful map
  \begin{equation*}\label{forgetful}
    \VZ^{\Gamma'}_{1,n+m}(\P^{\underline{n}},\partial \P^{\underline{n}})\to \VZ^\Gamma_{1,n}(\P^{\underline{n}},D),
  \end{equation*}
  there is an isomorphism
  $$\hat{\mathcal O}_{\VZ^{\Gamma}_{1,n}(\P^{\underline{n}},D),f}\xrightarrow{\cong} \hat{\mathcal O}_{\VZ^{\Gamma'}_{1,n+m}(\P^{\underline{n}},\partial \P^{\underline{n}}),f'}$$
   of complete local rings.\qed
\end{cor}

On the tropical side, the procedure of completing to toric boundary generalizes Example \ref{ex:genericonce}, by adding rays to the source curve while mapping them to new cones in the tropicalization $\Sigma_{(\P^{\underline n},\partial \P^{\underline n})}$. Any vertex of the source curve that is not in the $(n_1+\dots+n_a)$-dimensional cones corresponds to a component that is not contracted into the deepest intersection of boundary divisors. Therefore those components will react with the newly added divisors, by adding to their tropical vertices degree-many unbounded edges along each new direction of the missing boundary divisors $\partial\P^{\underline n}\setminus D$.

Once this is done, we obtain an honest balanced tropical map to a toric fan.

\begin{ex}\label{ex:ellipticcurve}
  (...continuing Example \ref{ex:genericonce}) We could further enhance the tropical map obtained in Example \ref{ex:genericonce} to one that maps to toric $\P^2$ by adding in the last missing divisor $D_2=\{x_2=0\}$.
  \begin{figure}[h]
    \centering
    \begin{tikzpicture}
      \tikzstyle{every node}=[font=\small]
      \draw [->] (-0.5,0)--(1.5,0);
      \draw [blue,->] (-0.5,0)--(1,0.1);
      \draw [blue,->] (-0.5,0)--(1,0);
      \draw [blue,->] (-0.5,0)--(1,-0.1);
      \draw [->] (2.5,0)--(4,0);
      \draw [->] (5,0)--(7,0);
      \draw [->] (5,0)--(5,2);
      \filldraw[blue] (-0.5,0) circle (1.5pt);
      \filldraw[blue] (5,0) circle (1.5pt);
      \draw [blue,->] (5,0)--(6.5,0.1);
      \draw [blue,->] (5,0)--(6.5,0);
      \draw [blue,->] (5,0)--(6.5,-0.1);
      \draw [blue,->] (5,0)--(5,1.5);
      \draw [blue,->] (5,0)--(5.1,1.5);
      \draw [blue,->] (5,0)--(4.9,1.5);
      \node at (-0.5,-0.4) [text=blue] {$g=1$};
      \node at (-0.5,-0.8) [text=blue] {$\text{degree}=3$};
      \node at (5,-0.4) [text=blue] {$g=1$};
      \node at (5,-0.8) [text=blue] {$\text{degree}=3$};
      \filldraw[blue] (11.5,0) circle (1.5pt);
      \draw [->] (8,0)--(9.5,0);
      \draw [->] (11.5,0)--(13.5,0);
      \draw [->] (11.5,0)--(11.5,2);
      \draw [->] (11.5,0)--(10.086,-1.414);
      \draw [blue,->] (11.5,0)--(13,0.1);
      \draw [blue,->] (11.5,0)--(13,0);
      \draw [blue,->] (11.5,0)--(13,-0.1);
      \draw [blue,->] (11.5,0)--(11.5,1.5);
      \draw [blue,->] (11.5,0)--(11.6,1.5);
      \draw [blue,->] (11.5,0)--(11.4,1.5);
      \draw [blue,->,rotate around={-135:(11.5,0)}] (11.5,0)--(13,0.1);
      \draw [blue,->,rotate around={-135:(11.5,0)}] (11.5,0)--(13,0);
      \draw [blue,->,rotate around={-135:(11.5,0)}] (11.5,0)--(13,-0.1);
      \node at (10.9,0) [text=blue] {$g=1$};
      \node at (-0.5,-2)[rotate=90] {$\leftrightsquigarrow$};
      \node at (-0.5,-3) {$f\in \VZ^{[1,1,1]}_{1,3}(\P^2,D_0)$};
      \node at (5.5,-2)[rotate=90] {$\leftrightsquigarrow$};
      \node at (5.5,-3) {$\tsup[1]f\in \VZ^{\Gamma'}_{1,6}(\P^2,D_0+D_1)$};
      \node at (11.5,-2)[rotate=90] {$\leftrightsquigarrow$};
      \node at (11.5,-3) {$\tsup[2]f\in \VZ^{\Gamma''}_{1,9}(\P^2,\partial \P^2)$};
    \end{tikzpicture}
    \caption{Lifting the tropical map to toric fan of $\Sigma_{\P^2}$.}
    \label{generic}
  \end{figure}
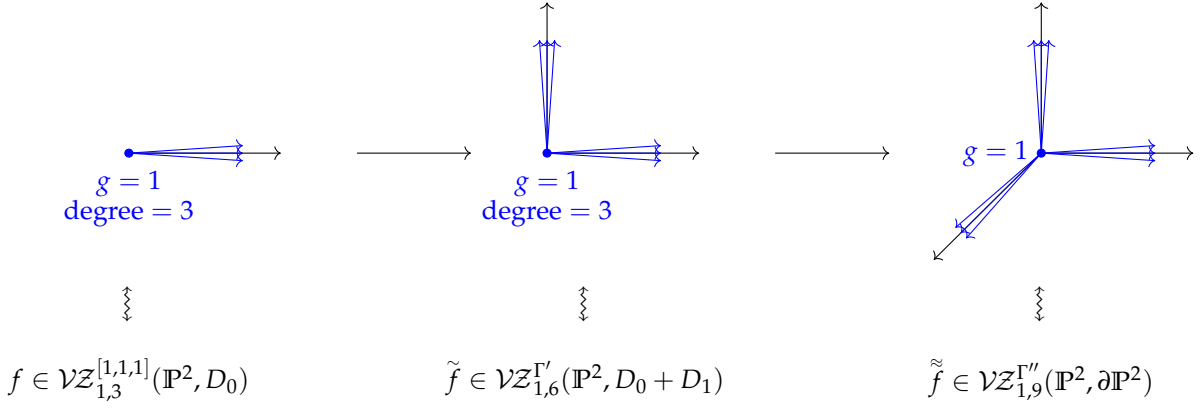
\end{ex}

Since the generic behavior for the tropicalization is to only add unbounded edges along specified directions, the tropical moduli for $f$, $\tsup[1]{f}$, $\tsup[2]{f}$ are exactly the same. Also because local neighborhood for $f$ and $\tsup[2]f$ are identical, analysis for $f$ can instead be done for $\tsup[2]{f}$. We use this idea to prove the logarithmic smoothness of $\VZ^\Gamma_{1,n}(\P^{\underline n},D)$ in the next section.

\subsection{Logarithmic smoothness}
As proved in the last section, we have a morphism $$\VZ^{\Gamma'}_{1,n+m}(\P^{\underline{n}},\partial \P^{\underline{n}})\to \VZ^\Gamma_{1,n}(\P^{\underline{n}},D)$$ satisfying Property ($\star$). Moreover, the space $\VZ^{\Gamma'}_{1,n+m}(\P^{\underline{n}},\partial \P^{\underline{n}})$ is logarithmically smooth by~\cite[Theorem 3.5.1]{RSW1}. Together we show $\VZ^\Gamma_{1,n}(\P^{\underline{n}},D)$ is also logarithmically smooth.

\begin{thm}\label{logsmooth}
  The relative genus one moduli spaces $\VZ^\Gamma_{1,n}(\P^{\underline{n}},D)$ are logarithmically smooth.
\end{thm}

\begin{pf}
  By Corollary \ref{cor:fulltoric}, there is a forgetful morphism
  $$\VZ^{\Gamma'}_{1,n+m}(\P^{\underline{n}},\partial \P^{\underline{n}})\to \VZ^\Gamma_{1,n}(\P^{\underline{n}},D)$$
  satisfying Property ($\star$), inducing an isomorphism
  $$\phi:\hat{\mathcal O}_{\VZ^{\Gamma}_{1,n}(\P^{\underline{n}},D),f}\xrightarrow{\cong} \hat{\mathcal O}_{\VZ^{\Gamma'}_{1,n+m}(\P^{\underline{n}},\partial \P^{\underline{n}}),f'}.$$
  By the analysis at the end of Section \ref{sec:fulltoric}, the map $\phi$ also sends the ideal of the boundary to the ideal of the boundary, because locally the tropical moduli remains unchanged.

  On the other hand, because $\VZ^{\Gamma'}_{1,n+m}(\P^{\underline{n}},\partial \P^{\underline{n}})$ is logarithmically smooth (i.e. toroidal), there exists, for $f'$, an affine toric variety and a point $p_{f'}\in W_{f'}$ with an isomorphism of complete local rings
  $$\psi:\hat{\mathcal O}_{\VZ^{\Gamma'}_{1,n+m}(\P^{\underline{n}},\partial \P^{\underline{n}}),f'}\cong \hat{\mathcal O}_{W_{f'},p_{f'}}$$ that sends the ideal of the boundary to that of the boundary $W_{f'}\setminus T$.

  Combining these two morphisms, we see for any $f\in\VZ^{\Gamma}_{1,n}(\P^{\underline{n}},D)$, there is an isomorphism of complete local rings
  $$\psi\circ \phi:\hat{\mathcal O}_{\VZ^{\Gamma}_{1,n}(\P^{\underline{n}},D),f}\xrightarrow{\cong}\hat{\mathcal O}_{W_{f'},p_{f'}}$$
  that sends the ideal of the boundary to the boundary. Therefore $\VZ^{\Gamma}_{1,n}(\P^{\underline{n}},D)$ is also logarithmically smooth.\qed
\end{pf}

\begin{prop}
  The moduli space $\VZ^{\Gamma}_{1,n}(\P^{\underline{n}},D)$ carries a well-defined virtual fundamental class in the expected degree.
\end{prop}
\begin{pf}
  Because the space $\VZ^{\Gamma}_{1,n}(\P^{\underline{n}},D)$ is logarithmically smooth, it is irreducible and is of expected dimension. Therefore taking its fundamental class suffices.\qed
\end{pf}

We will use the fundamental class of $\VZ^{\Gamma}_{1,n}(\P^{\underline{n}},D)$ to define a virtual class for the space of genus 1 maps to a more general target pair in the last section.

\subsection{An example comparison between the genus one moduli spaces}\label{sec:vzexample}
We devote this short section to a calculation to show the difference between the factorization property and the well-spacedness property.

We have three moduli spaces here:
\begin{itemize}
  \item the space $\widetilde{\VZ}_{1,n}(\P^{\underline n},\boldsymbol{d})$ of logarithmic stable maps from radial genus 1 and $n$-pointed logarithmic curves to $\P^{\underline n}$ of degree $\boldsymbol{d}$;
  \item the space $\VZ_{1,n}(\P^{\underline n},\boldsymbol{d})$ of maps in $\widetilde{\VZ}_{1,n}(\P^{\underline n},\boldsymbol{d})$ that satisfy the factorization property; and
  \item the space $\VZ_{1,n}(\P^{\underline n},\emptyset)$ of maps in $\widetilde{\VZ}_{1,n}(\P^{\underline n},\boldsymbol{d})$ that satisfy the well-spacedness condition.
\end{itemize}
The space $\VZ_{1,n}(\P^{\underline n},\boldsymbol{d})$ is characterized by the factorization property defined and investigated in \cite{RSW1}. When the target is a projective space, it was shown to resemble the original Vakil--Zinger desingularization of the moduli space of genus one stable maps. However, if the target is a product of projective spaces $\P^{\underline n}$, the space $\VZ_{1,n}(\P^{\underline n},\boldsymbol{d})$ is not smooth.

The space $\VZ_{1,n}(\P^{\underline n},\emptyset)$ is the one constructed earlier in this paper. By Theorem \ref{logsmooth}, this space is logarithmically smooth.

\begin{ex}
  Consider the space of 0-marked genus 1 curves to $\P^1\times \P^1$ with bidegree $(2,2)$. By the degree-genus formula, a smooth projective curve of bidegree $(2,2)$ in $\P^1\times \P^1$ will have genus $(2-1)\cdot(2-1)=1$.

  By Riemann--Roch theorem, the expected dimension of this moduli space should be:
  \begin{align*}
    &\int_\beta c_1(\P^1\times \P^1)+(\dim(\P^1\times \P^1) - 3)(1-g)+n\\
    =& 2\cdot 2+2\cdot 2+0+0= 8.
  \end{align*}
  So the main component of maps from smooth curves will have dimension 8.

  Consider another stratum: the source curve consists of a genus 1 component and a genus 0 component. The map will contract the genus 1 component into a fiber of one factor of $\P^1$, and contract the genus 0 component to a fiber of the other factor of $\P^1$ as a double cover. See Figure \ref{moduliexample} for an illustration.

  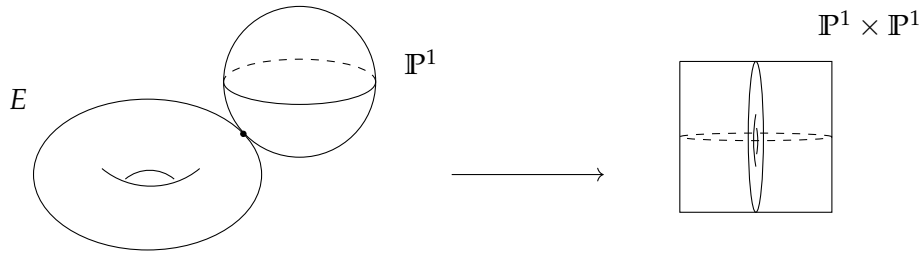
\begin{figure}[h]
    \centering
    \captionsetup{justification=centering,margin=2cm}
      \begin{tikzpicture}[scale=1]
        \draw (0,0) ellipse (1.5cm and 1cm);
        \draw (-0.6,0.08) arc (-130:-50:1);
        \draw (0.35,-0.06) arc (50:130:0.5);
        \draw (2,1.23) circle (1);
        \draw[dashed] (3,1.23) arc (0:180:1 and 0.3);
        \draw (3,1.23) arc (0:-180:1 and 0.3);
        \filldraw (1.26,0.54) circle (1pt);
        \node at (-1.7,1) {$E$};
        \node at (3.6,1.5) {$\P^1$};
        \draw [->]  (4,0) -- (6,0);
        \node at (9.5,2) {$\P^1\times \P^1$};
        \draw (7,1.5) -- (9,1.5) -- (9,-0.5) -- (7,-0.5)--(7,1.5);
        \draw (8,0.5) ellipse (0.1cm and 1cm);
        \draw (8,0.8) arc (-190:-170:2);
        \draw (8.01,0.27) arc (-10:10:1);
        \draw [dashed] (8,0.5) ellipse (1cm and 0.05cm);
        \end{tikzpicture}
        \caption{The stratum of maps of the required shape.}
      \label{moduliexample}
  \end{figure}
  Here the dashed line represents the double cover of a fiber of $\P^1$ along the vertical direction. We calculate the dimension of this stratum. We start by choosing a point on $\P^1$, the fiber of which will be the target of the genus 1 component. Upon fixing the point, the moduli for the genus 1 component is then $M_{1,1}(\P^1,2)$, which is of dimension 5. The genus 0 component's contribution is essentially $M_{0,1}(\P^1\times \P^1,(0,2))$, which is of dimension 4.

  Finally we are gluing the two marked points, which corresponds to taking the fiber product over $\P^1\times \P^1$ (where maps to it are the evaluation maps on the marked points). Together we have the dimension of this stratum:
  $$\dim S = (5+1) + 4 - 2=8.$$

  Notice that this stratum lives in $\VZ_{1,n}(\P^{\underline n},\boldsymbol{d})$, as it does not contract the genus 1 component entirely as a map to $\P^1\times \P^1$. As a result, there are two 8-dimensional pieces in this moduli space $\VZ_{1,n}(\P^{\underline n},\boldsymbol{d})$, neither a degeneration of the other. In particular, the space $\VZ_{1,n}(\P^{\underline n},\boldsymbol{d})$ is not logarithmically smooth. So the well-spacedness condition is a non-trivial requirement for the desingularization.
\end{ex}

\section{Very ample snc pairs of the same degree and the virtual pullbacks}
This section has two parts. The first part constructs a virtual class for the space $\VZ^\Gamma_{1,n}(X,Y)$, where $(X,Y)$ is a very ample SNC pair of the same degree, by mapping to projective spaces and also the technique of virtual pullbacks. This generalizes the construction in~\cite{BNR} where $Y$ is a smooth divisor.

The second part proves the fact that the map forgetting fictitious component and markings is birational, therefore relating the virtual classes. This means we could iteratively compute the virtual classes for certain genus 1 moduli spaces entirely within the scope of the genus 1 theory built in this paper.

\subsection{Virtual class for a very ample SNC pair of the same degree}\label{sec:generalpair}
Let $(X,Y)$ be an SNC pair of the same degree, meaning that $Y\subset X$ is a simple normal crossings divisor and also all components of $Y$ are defined by sections of the same line bundle.

The moduli space $\VZ^\Gamma_{1,n}(X,Y)$ is the space of maps from radially aligned genus 1 curves to $(X,Y)$ with contact orders $\Gamma$ that also satisfy the well-spacedness condition, constructed analogously to $\VZ^\Gamma_{1,n}(\P^{\underline n},D)$ explained in Section \ref{sec:moduli}.

Because of the assumption on $(X,Y)$, there is an embedding $X\to \P^{r}$ and thus a map $$\pi:\VZ^\Gamma_{1,n}(X,Y) \to \VZ^\Gamma_{1,n}(\P^{r},D=\sum H_i).$$ Let $\VZ_{1,n}(X,\beta)$ be the moduli space of radial maps to $X$ that only satisfy the factorization condition (Definition \ref{def:faccon}). There is a forgetful map
$$\VZ^\Gamma_{1,n}(X,Y)\to \VZ_{1,n}(X,\beta)$$
obtained by forgetting all relative information. Similarly there is another forgetful map
$$s:\VZ^\Gamma_{1,n}(\P^{r},D)\to \VZ_{1,n}(\P^{r},\boldsymbol{d}).$$

Together, we have the following commutative diagram, which is also Cartesian in both the category of fine and saturated logarithmic stacks and coherent logarithmic stacks (because $t$ is strict, as the logarithmic structures are both pulled back from $\mathfrak M^{\rad}_{1,n}$):
\begin{center}
  \begin{tikzcd}
    \VZ^\Gamma_{1,n}(X,Y) \arrow[r,"\pi"]\arrow[d] \arrow[dr, phantom, "\square"] & \VZ^\Gamma_{1,n}(\P^{r},D) \arrow[d,"s"]\\
    \VZ_{1,n}(X,\beta) \arrow[r,"t"] &\VZ_{1,n}(\P^{r},\boldsymbol{d}).
  \end{tikzcd}
\end{center}

The above Cartesian square is equivalent to another Cartesian square
\begin{center}
  \begin{tikzcd}
    \VZ^\Gamma_{1,n}(X,Y) \arrow[r]\arrow[d] \arrow[dr, phantom, "\square"]&\VZ^\Gamma_{1,n}(\P^{r},D)\times \VZ_{1,n}(X,\beta)\arrow[d,"s\times t"]\\
     \VZ_{1,n}(\P^{r},\boldsymbol{d})\arrow[r,"\Delta"]& \VZ_{1,n}(\P^{r},\boldsymbol{d})\times \VZ_{1,n}(\P^{r},\boldsymbol{d}),
  \end{tikzcd}
\end{center}
where $\Delta$ is the diagonal map. Since the space $\VZ_{1,n}(\P^{r},\boldsymbol{d})$ is smooth by Theorem \ref{thm:abssmooth}, the diagonal map $\Delta$ is a regular embedding. By~\cite[Example 8.3.1]{Fulton}, there is a well-defined pullback map
$$\Delta^!:A_*(\VZ^\Gamma_{1,n}(\P^{r},D)\times \VZ_{1,n}(X,\beta))\to A_*(\VZ^\Gamma_{1,n}(X,Y)).$$

The space $\VZ_{1,n}(X,\beta)$ has a virtual class in expected degree by~\cite[Theorem 4.4.1]{RSW1}, which can be used to define the virtual class on $\VZ^\Gamma_{1,n}(X,Y)$ by pulling back its product with $[\VZ^\Gamma_{1,n}(\P^{r},D)]$.
\begin{thm}
  The moduli space $\VZ^\Gamma_{1,n}(X,Y)$ carries a well-defined virtual class in expected degree, defined by
  $$[\VZ^\Gamma_{1,n}(X,Y)]^\vir:=\Delta^!([\VZ^\Gamma_{1,n}(\P^{r},D)]\times [\VZ_{1,n}(X,\beta)]^\vir).$$
\end{thm}

\begin{rem}
  One would ideally want to replace $\VZ_{1,n}(X,\beta)$ with $\VZ_{1,n}(X,\emptyset)$ in the Cartesian diagram and therefore show the results for any very ample SNC pair without the same degree restraints. However, it is not immediately possible here because we do not have a virtual class on the space $\VZ_{1,n}(X,\emptyset)$.
\end{rem}
\subsection{Birationality when removing fictitious component and markings}\label{forget}
Let $\VZ^\Gamma_{1,n}(\P^{\underline{n}},D)$ be the moduli space of maps to $(\P^{\underline{n}},D)$. Write $\{p_1,\dots,p_n\}$ for the marked points and $D=\sum_{i=1}^s D_i$ for the irreducible components of $D$. The number $\Gamma_{ij}$ is the contact order of $p_i$ with $D_j$.

Suppose there is a number $m<n$ such that the contact orders in $\Gamma$ satisfy the following:
\begin{itemize}
  \item $\Gamma_{i1}=1$ for $1\le i\le m$ and $\Gamma_{i1}=0$ for $i>m$; namely only $\{p_1,\dots,p_m\}$ intersect $D_1$, and they intersect transversely; and
  \item $\Gamma_{ij}=0$ for $1\le i\le m$ and $j>1$; namely those $m$ points do not intersect any other $D_j$.
\end{itemize}
In this case, we call the $m$ marked points together with the component $D_1$ \textit{fictitious}, as we will prove shortly that they will not meaningfully contribute towards the virtual class of the space.

There exists a composition of forgetful maps:
\begin{center}
  \begin{tikzcd}
    \VZ^\Gamma_{1,n}(\P^{\underline{n}},D) \arrow["f"]{r}\arrow[rr,bend left=25,"h"]{rr} & \VZ^{\Gamma''}_{1,n}(\P^{\underline{n}},D-D_1) \arrow["g"]{r} &\VZ^{\Gamma'}_{1,n-m}(\P^{\underline{n}},D-D_1),
  \end{tikzcd}
\end{center}
where $f$ forgets all relative information with respect to $D_1$ and $g$ forgets the marked points $p_i$ for $1\le i \le m$. The symbols $\Gamma''$ and $\Gamma'$ represent the respective submatrices of $\Gamma$.

\begin{lem}\label{birational}
  The map $h:\VZ^\Gamma_{1,n}(\P^{\underline{n}},D)\to \VZ^{\Gamma'}_{1,n-m}(\P^{\underline{n}},D-D_1)$ relates the virtual fundamental classes by the identity
  $$h_*[\VZ^\Gamma_{1,n}(\P^{\underline{n}},D)]=m!\cdot[\VZ^{\Gamma'}_{1,n-m}(\P^{\underline{n}},D-D_1)].$$
\end{lem}
\begin{pf}
  We have shown that both spaces are logarithmically smooth by Theorem \ref{logsmooth}. So it suffices to check $h$ is $m!:1$ on the locus of maps with trivial logarithmic structures. This is immediate because $m!$ is the number of permutations of the $m$ markings.\qed
\end{pf}

\begin{ex}
  We consider genus 0 stable maps for simplicity. Let $\overline{M}_{0,2}^{(1,1)}(\P^2\mid L,2)$ be the space of 2-marked, genus 0, degree 2 stable maps to $\P^2$ relative to a line $L$ with transverse contacts. There is a composition of forgetful maps to the spaces of Kontsevich stable maps
  \begin{center}
    \begin{tikzcd}
      \overline{M}_{0,2}^{(1,1)}(\P^2\mid L,2) \arrow["f"]{r}\arrow[rr,bend left=25,"h"]{rr} & \overline{M}_{0,2}(\P^2,2) \arrow["g"]{r} &\overline{M}_{0,0}(\P^2,2).
    \end{tikzcd}
  \end{center}
For a quick dimension check, $\dim \overline{M}_{0,2}^{(1,1)}(\P^2\mid L,2)=5$ by realizing the open interior as generic conics (thus not tangent to $L$). Also
$$\dim \overline{M}_{0,2}(\P^2,2)=2+c_1(\P^2)\cdot 2H-1=7,$$
and by forgetting the 2 marked points, $\dim \overline{M}_{0,0}(\P^2,2)=5=\dim \overline{M}_{0,2}^{(1,1)}(\P^2\mid L,2)$.

For an element in the open locus of $\overline{M}_{0,0}(\P^2,2)$, we expect it to meet $L$ in two distinct points. So generically $h$ is $2:1$.
\end{ex}

It is possible to lift the above calculations to a very ample SNC pair $(X,Y)$ of the same degree. Suppose $Y_1\subset Y$ and $\{p_1,\dots,p_m\}$ are fictitious, as defined in the beginning of this section. There exists a Cartesian square
\begin{equation}
  \begin{tikzcd}
    \VZ^\Gamma_{1,n}(X, Y) \arrow["p"]{r} \arrow["g"]{d} \arrow[dr, phantom, "\square"]& \VZ^\Gamma_{1,n}(\P^{r}, D)\arrow["h"]{d}\\
    \VZ^{\Gamma'}_{1,n-m}(X, Y-Y_1) \arrow["q"]{r} & \VZ^{\Gamma'}_{1,n-m}(\P^{r}, D-D_1),
  \end{tikzcd}
  \tag{$\dagger$}
\end{equation}
where $\Gamma'$ is the submatrix of contact orders of $\{p_{m+1},\dots,p_n\}$ with $Y-Y_1$. More precisely, this happens iff
\begin{itemize}
  \item the only points that intersect $Y_1$ are $\{p_1,\dots,p_m\}$, and they intersect transversely; and
  \item $\{p_1,\dots,p_m\}$ do not intersect any other component of $Y$.
\end{itemize}

\begin{thm}
  For a very ample SNC pair $(X, Y)$ of the same degree, if $Y_1$ and $\{p_1,\dots,p_m\}$ are fictitious inducing the diagram $(\dagger)$, then the forgetful map
  $$g:\VZ^\Gamma_{1,n}(X, Y)\to \VZ^{\Gamma'}_{1,n-m}(X, Y-Y_1)$$
  relates the virtual classes by the following identity:
  $$g_*[\VZ^\Gamma_{1,n}(X, Y)]^{\text{vir}}=m!\cdot [\VZ^{\Gamma'}_{1,n-m}(X, Y-Y_1)]^{\text vir}.$$
\end{thm}

\begin{note}
  It is worth noting that these results are not obvious without some positivity conditions on the divisor $Y$; there are some interesting examples discussed by Tehrani and Zinger in ~\cite{Teh}.
\end{note}

\begin{pf}
  For simplicity of notations, we prove this result for the virtual classes defined on $\VZ_{1,n}(\P^{r}, \boldsymbol{d})$. The result for refined classes defined on $\VZ^\Gamma_{1,n}(X, Y)$ follows similarly.

  There are two Cartesian diagrams
  \begin{center}
    \begin{tikzcd}
      \VZ^\Gamma_{1,n}(X,Y) \arrow[r]\arrow[d] \arrow[dr, phantom, "\square"]&\VZ^\Gamma_{1,n}(\P^{r},D)\arrow[d]& \VZ^{\Gamma'}_{1,n-m}(X,Y-Y_1) \arrow[r]\arrow[d] \arrow[dr, phantom, "\square"]&\VZ^{\Gamma'}_{1,n-m}(\P^{r},D-D_1)\arrow[d]\\
      \VZ_{1,n}(X,\beta) \arrow[r] &\VZ_{1,n}(\P^{r},\boldsymbol{d}),&\VZ_{1,n-m}(X,\beta) \arrow[r] &\VZ_{1,n-m}(\P^{r},\boldsymbol{d}).
    \end{tikzcd}
  \end{center}
  We assume that
  \begin{alignat*}{3}
    &[\VZ^\Gamma_{1,n}(X,Y)]^\vir &&= [\VZ_{1,n}(X,\beta)]^\vir\cap [\VZ^\Gamma_{1,n}(\P^{r}, D)]&&\in A_*(\VZ_{1,n}(\P^{r},\boldsymbol{d})),\\
    &[\VZ_{n-m}^{\Gamma'}(X, Y-Y_1)]^\vir &&= [\VZ_{1,n-m}(X,\beta)]^\vir\cap [\VZ^{\Gamma'}_{1,n-m}(\P^{r},D-D_1)]&&\in A_*(\VZ_{1,n-m}(\P^{r},\boldsymbol{d})).
  \end{alignat*}

  Then it suffices to show the forgetful map
    $$g':\VZ_{1,n}(\P^{r}, \boldsymbol{d})\to \VZ_{1,n-m}(\P^{r},\boldsymbol{d})$$
  pushes $[\VZ^\Gamma_{1,n}(X,Y)]^\vir$ to a multiple of $[\VZ_{1,n-m}^{\Gamma'}(X, Y-Y_1)]^\vir$.

  Note that since $g'$ is flat, the pullback of $[\VZ_{1,n-m}(X,\beta)]^\vir$ along $g'$ is $[\VZ_{1,n}(X,\beta)]^\vir$. Then using the projection formula and Lemma \ref{birational},
  \begin{align*}
    g'_*([\VZ^\Gamma_{1,n}(X,Y)]^\vir) &=g'_*([\VZ_{1,n}(X,\beta)]^\vir\cap [\VZ^\Gamma_{1,n}(\P^{r}, D)])\\
    &= g'_*(g'^*[\VZ_{1,n-m}(X,\beta)]^\vir\cap [\VZ^\Gamma_{1,n}(\P^{r}, D)])\\
    &=[\VZ_{1,n-m}(X,\beta)]^\vir\cap g'_*([\VZ^\Gamma_{1,n}(\P^{r}, D)])\\
    &=[\VZ_{1,n-m}(X,\beta)]^\vir\cap m!\cdot[\VZ^{\Gamma'}_{1,n-m}(\P^{r},D-D_1)]\\
    &=m!\cdot[\VZ_{n-m}^{\Gamma'}(X, Y-Y_1)]^\vir.
  \end{align*}
  Therefore, we have the required identity.\qed
\end{pf}

\bibliographystyle{alpha}
\bibliography{ref.bib}

\end{document}